**Aleksandr RAKHMANGULOV\*, Anatoly KOLGA, Nikita OSINTSEV**
Nosov Magnitogorsk State Technical University
Pr. Lenina, 38, 455000, Magnitogorsk, Russia
**Ivan STOLPOVSKIKH**
Kazah National Technical University named after K.I. Satpayev,
Satpaev 22, 050013, Almaty, Kazakhstan
**Aleksander SŁADKOWSKI**
Silesian University of Technology, Faculty of Transport
Krasinski 8, 40-019 Katowice, Poland
*Corresponding author*. E-mail: prtrans@gmail.com

# MATHEMATICAL MODEL OF OPTIMAL EMPTY RAIL CAR DISTRIBUTION AT RAILWAY TRANSPORT NODES

**Summary.** At present there are two trends in the market of rail freight transportation in Russia: freight owners put forward higher demands to the transportation quality (promptness of delivery) in an effort to reduce storage costs by means of reducing the size of freight shipment; the structure of railcar traffic volume of the railways of Russia is getting more complex due to the reduction of the average shipment size and due to the transfer of railcar fleet ownership to a large number of operating companies. These trends significantly complicate operational management supervision of railway stations and transport nodes. Application of typical data from the information system about the railcar location at the transportation node is not enough for the dispatchers to make the best decision concerning the car traffic management. The dispatcher traffic control service needs some software-based models of efficient railcar distribution. The article is concerned with the description and development of the mathematical model of empty railcar distribution for loading at the railway transport node; this model will take into account the requirements of railcar owners in terms of their cars application, the operating work level of railroad stations of the transportation node and the possibility of adding the groups of empty railcars to the transfer trains, clean-up trains and industrial railway trains operating on a tight schedule. The developed model and the software package were implemented in the information system of the industrial railway of the major metallurgical enterprise - OJSC «Magnitogorsk Metallurgical Works», which processes up to two thousand of railcars belonging to different owners. This model made it possible to reduce the labour intensity of dispatcher operation planning the empty railcar distribution for loading and reduce the total time the railcars spend in the enterprise railway system.

# МАТЕМАТИЧЕСКАЯ МОДЕЛЬ ОПТИМАЛЬНОГО РАСПРЕДЕЛЕНИЯ ПОРОЖНИХ ВАГОНОВ В ЖЕЛЕЗНОДОРОЖНЫХ ТРАНСПОРТНЫХ УЗЛАХ

**Аннотация.** Двумя основными тенденциями рынка железнодорожных перевозок в России являются: повышение требований грузовладельцев к качеству (своевременности) перевозок, что, в частности, связано с их стремлением к сокращению складских издержек путем уменьшения размера транспортно-грузовых партий; усложнение структуры вагонопотоков на железных дорогах России, происходящее как по причине уменьшения среднего размера отправок, так и в результате передачи парка железнодорожных вагонов в собственность множеству операторских компаний. Отмеченные тенденции значительно усложняют оперативное руководство работой железнодорожных станций и транспортных узлов. Использование диспетчерами типовых данных информационной системы о местоположении вагонов в транспортном узле является недостаточным для выработки оптимального решения по управлению вагонопотоками. Диспетчерскому аппарату требуется программно реализованные модели оптимального использования вагонов. В статье рассматривается постановка и алгоритм реализации математической модели распределения порожних вагонов под погрузку в железнодорожном транспортном узле, учитывающей требования собственников вагонов на их использование, оперативный уровень загрузки железнодорожных станций узла и возможность включения групп порожних вагонов в состав передаточных, вывозных поездов и поездов, обращающихся по контактному графику. Разработанная модель и компьютерная программа реализованы в рамках информационной системы железнодорожного транспорта крупнейшего в мире металлургического предприятия – ОАО «Магнитогорский металлургический комбинат», ежесуточно перерабатывающего до двух тысяч вагонов, принадлежащих различным собственникам. Использование модели позволило значительно сократить трудоемкость оперативного планирования работы диспетчеров по распределению порожних вагонов под погрузку, сократить суммарное время нахождения вагонов на путях предприятия.

## 1. INTRODUCTION

Reorganization of the transport system and, particularly, the federal railway (in Russia) restructuring caused significant changes in the structure of the transportation services, such as: the increase in the number of haulage companies including carriers, which have their own rolling stock. The main trend of the transportation market is that customers put forward higher demands in terms of the transport service quality. To reduce costs customers try to reduce the consignment size and improve the promptness of delivery.

The trends above are the objective causes of freight traffic flow sophistication on the one hand, and the increase in the requirements to the quality of traffic control, on the other hand [1].

One of the main factors, which makes it difficult to improve the quality of freight traffic at the time of freight traffic flow sophistication, is the lack of co-ordination in the work of the operating managers of the transportation process. Large metallurgical enterprises have developed railway systems and a large amount of rolling stock. For this kind of transport are solving complicated logistical problems [2]. Only inefficient interaction between the main line and the industrial railway transport of metallurgical enterprises results in annual logistical costs of up to 1.5 billion rubles ($45 mln. dollars) on average per enterprise [3, 4]. One of the pressing problems caused by the transfer of railcars into haulage companies' ownership is the rational use of the empty railcar fleet at transportation nodes taking into account various restrictions, which railcar owners face in terms of the cars application [5].

It is necessary to be able to solve this problem on-line as the data concerning the allowable directions and consignments for railcars of various owners can vary daily. Practical experience shows that empty railcar distribution for loading takes place at the beginning of the design day and, as a rule, it is not adjusted until the beginning of the next base period (schedule day). This feature makes it possible to view the problem of optimal distribution of empty railcars as static one and to simplify it to an allocation linear programming problem, in particular, to a transport problem.

On the other hand, account of additional restrictions for the owners on the use of empty railcars connected with the loading operations of certain cargos and shipment of railcars to certain destinations does not differ much from the well known allocation problems for various railcars over the loading areas [6]. That is why the authors believe that this problem should be solved taking into account that the groups of empty railcars are bound to the schedules of district and transfer trains operation or to the tight schedule of trains operating between the stations of the enterprise. The restrictions above make it possible to solve this problem in the on-line mode only on the basis of the developed mathematical model, which uses the data from the information system to distribute the empty railcars on-line [7 - 9].

## 2. STATEMENT OF THE PROBLEM OF EMPTY RAILCAR DISTRIBUTION AT RAILROAD TRANSPORTATION NODES

Statement of the problem of optimal distribution of empty railcars at railroad transportation nodes as a static transport problem can be formulated as follows. By the beginning of the design day there are $k = 1, 2, \ldots, L$ groups of empty railcars of various kinds and belonging to different owners on the on the special railway track. Here we will denote the number of railcars belonging to each group as $A_k$.

By the start of the base period the railcars belonging to each group $k$ can be located at different industrial railway stations. We will assign to each industrial station the index $i = 1, 2, \ldots, M$, where $M$ is the number of stations making up the railway transport node. We will denote a part of railcars from the $k$ group located at the $i$ station as $A_{ki}$. It is necessary to distribute the empty railcars belonging to different groups among the industrial stations.

We will denote the destination stations (loading stations) as $j = 1, 2, \ldots, N$ indices, where $N$ is the number of empty railcar consumers. The number of railcars belonging to $k$ group, which must be delivered to the corresponding stations will be denoted as $B_{kj}$.

The plan of empty railcar distribution must meet the requirement of complete distribution of the empty stock located on the approach lines, i.e.

$$\sum_{i=1}^{M} A_{ki} = \sum_{j=1}^{N} B_{kj} = A_k, k = 1, 2, \ldots, L. \tag{1}$$

In order to simulate the return of excess cars or the order of lacking empty railroad cars, the extra supplier and consumer corresponding to the connecting station were introduced.

It is necessary to determine the number of empty railcars belonging to each group $k$ and located at the stations $i$, which must be delivered for loading to points $j$. We will denote the unknown models, which as a whole form the plan of empty railcar distribution, as $x_{kij}$.

If the required moments of spotting are not predetermined, it is reasonable to select the minimum total time spent on delivery of empty railcars from the stations where they are located at the beginning of the base period to the loading stations as the optimality criterion. We will denote the time spent on each spotting $x_{kij}$ as $C_{kij}$.

The time consumption $C_{kij}$ is calculated using well known methods of the shortest path search for the transport network as a sum of estimations $p_{ij}$ of transport network paths forming the traffic route of the railroad cars from the starting station of the route $i$, to the final station $j$. Estimations $p_{ij}$ of the transport network paths (span) connecting the vertexes of the network (railway stations) can be calculated by the following formula

$$p_{ij} = t_i \sigma_i + t_{ij}, \tag{2}$$

where $t_i$ is the processing time (dead time) of empty railroad cars on the $i$ station (the starting station of the route). In the model it is set as the average time the transit railcars spend at the station;

$\sigma_i$ is the station workload factor, the variable that depends on the number of cars present at the station at the design moment and on the number of shunting facilities available. The value of the station workload factor can be determined using fuzzy logic methods described in the articles [4, 10];

$t_{ij}$ is the specified travelling time along the span connecting the starting station $i$ and the final station $j$ of the route.

Then, if the traffic route for empty railcars is determined $S_{ij} = \{i, ..., \lambda_i, j\}$, where $\lambda_j$ is the number of the vertex preceding the $j$-th one, time consumption $C_{kij}$ is calculated as the sum of the curve potentials $p_{ij}$ comprising the traffic route $S_{ij}$ of empty railcars $x_{kij}$.

The target function of railcar distribution will take the following form

$$\sum_{k=1}^{L}\sum_{i=1}^{M}\sum_{j=1}^{N} C_{kij} \cdot x_{kij} \to \min. \tag{3}$$

The target function (3) combined with the limitation (1) and the nonnegativity constraint of the number of car supply $x_{kij}$ form a known static multicommodity transport problem.

However, application of this simplified model for improving empty railcar distribution will result in a non-optimal plan, as it does not take into account the restrictions imposed by the train schedule (including tight train schedule of industrial railways). If we take into account the addition of empty railcars to the formation of trains operating according to the tight schedule, the period of time, which the empty railcars spend at stations on their route as well as the total time consumption $C_{kij}$, will increase due to extra waiting time of the scheduled train departure. Let us assume that for each $i$-th station the scheduled train departure time $t_{ir}$ is known, where $r = 1,2,...,R$ is the progressive number of the train leaving station $i$ for station $j$, $R$ is the number of trains outgoing from station $j$ over the base period of time. Then the formula (2) will take the following form

$$p_{ij} = (t_{ir} - p_i) + t_{ij}, \tag{4}$$

where $p_i$ is the assessment (potential) of the $i$-th vertex of the transport network. Variable $p_i$ is the total travel time from the starting station of the route to the $i$-th station and it is calculated as the sum of assessments (potentials) of the transport network curves making up the traffic route form the starting vertex of the route to the $i$-th one according to the following formula

$$p_i = p_{\lambda_i} + p_{ij}, \tag{5}$$

the difference $t_{ir} - p_i$ determines the dead time of the railcar delivery in expectation of the dispatch of the next scheduled train, which is supposed to include the empty railcar supply.

It is necessary to choose such a value of $r$ in the formula (4), which would satisfy the following conditions

$$t_{ir-1} < (p_i + t_i \sigma_i) \le t_{ir} \text{ и} \tag{6}$$

$$x_{kij} \le Q_{ir}, \tag{7}$$

where $Q_{ir}$ is the maximum number of railcars, which might be included into the train $r$ dispatched according to the schedule from the station $i$ at the time point $t_{ir}$.

Fulfilling of condition (6) makes sure that the dispatch of the empty railcars from the station will not take place before the dispatch of the nearest scheduled train and fulfilling of condition (7) will ensure the train size limitation.

The introduction of additional limitations (6) and (7) into the model in most cases makes it impossible to solve this problem for one iteration as in some situations it is impossible to provide delivery of the optimum number of empty railcars as part of the first trains at the beginning of the base period. This situation may arise when the size of railcar supply exceeds the value of $Q_{ir}$ or even the train size. In this case the remaining part of empty railcars, which will further be referred to as

«undistributed railcar supply», should be included into trains scheduled for later dispatch. This, according to the formula (4), will make it necessary to change the potentials of the transport network curves and the time consumption $C_{kij}$ and to make a new plan of empty railcar distribution.

Thus, in order to find the minimum of the target function (3) with limitations (1) and (4) it is necessary to develop an algorithm that will provide both optimum railcar distribution over the loading areas and the loading stations and optimum empty railcar assignment to the trains operating according to the tight schedule. The developed algorithm is given in Fig. 1 and it consists of the following operations:

1. Initial data preparation including the values of the following variables:

$A_k$ – the number of empty railcars of each group $k$ located at the railway transport node at the beginning of the base period, railcars;

$A_{ki}$ – the number of empty railcars of each group $k$ located at each industrial railway station $i$, railcars;

$B_{kj}$ – demand for empty railcars at each $j$-th station or loading area, railcars;

$t_i$ – the average processing time of transit rail cars at the $i$-th station, min.;

$\sigma_i$ – station workload factor. It is calculated according to a special algorithm on the basis of the data about the number of railcars at each station at the start of the base period (the data can be requested from the information system), the number of shunting locomotives and draw out tracks;

$t_{ij}$ – train travel time between the neighbouring industrial railway stations, min.;

$t_{ir}$ – train schedule – train dispatch time for each $r$-th railway station, min.;

$Q_{ir}$ – the maximum number of railcars, which can be included into the train $r$ dispatched according to the tight schedule from station $i$ at the time $t_{ir}$, railcars.

2. Calculation of the time spent on railcar processing at the stations $t_i \sigma_i$ for further use in the formula (6).

3. Development of optimal route set for all starting vertexes of the transport network (carried out according to a special algorithm). The $i$-th stations where railcars $A_{ki}$ are located are taken as the starting vertexes of the transport network. In the process of the optimal route set development it is necessary to adjust the curve assessments according to the formula (4) taking into account the limitations (6). Limitation (7) is not taken into account on this stage as the size of empty railcar supply $x_{kij}$ is not known yet.

4. Calculation of time $C_{kij}$, which is necessary to deliver the empty railcars from the starting station of each route $i$ to the final $j$-th stations where there is a demand for empty railcars $B_{kj}$. The value of $C_{kij}$ will be equal to the final vertex potential of the corresponding route, i.e. $C_{kij} = p_j$.

5. Successive solution of the static transport problem of linear programming as a matrix problem (1), (2) for each group of empty railcars $A_k$ and the calculation of the optimal railcar supply $x_{kij}$.

6. Route classification by the value of $C_{kij}$ in the ascending order and sequential condition testing (7) for all vertexes of each route. The vertexes are tested successively beginning with the $i$-th one to the $\lambda_j$-th one preceding the final vertex.

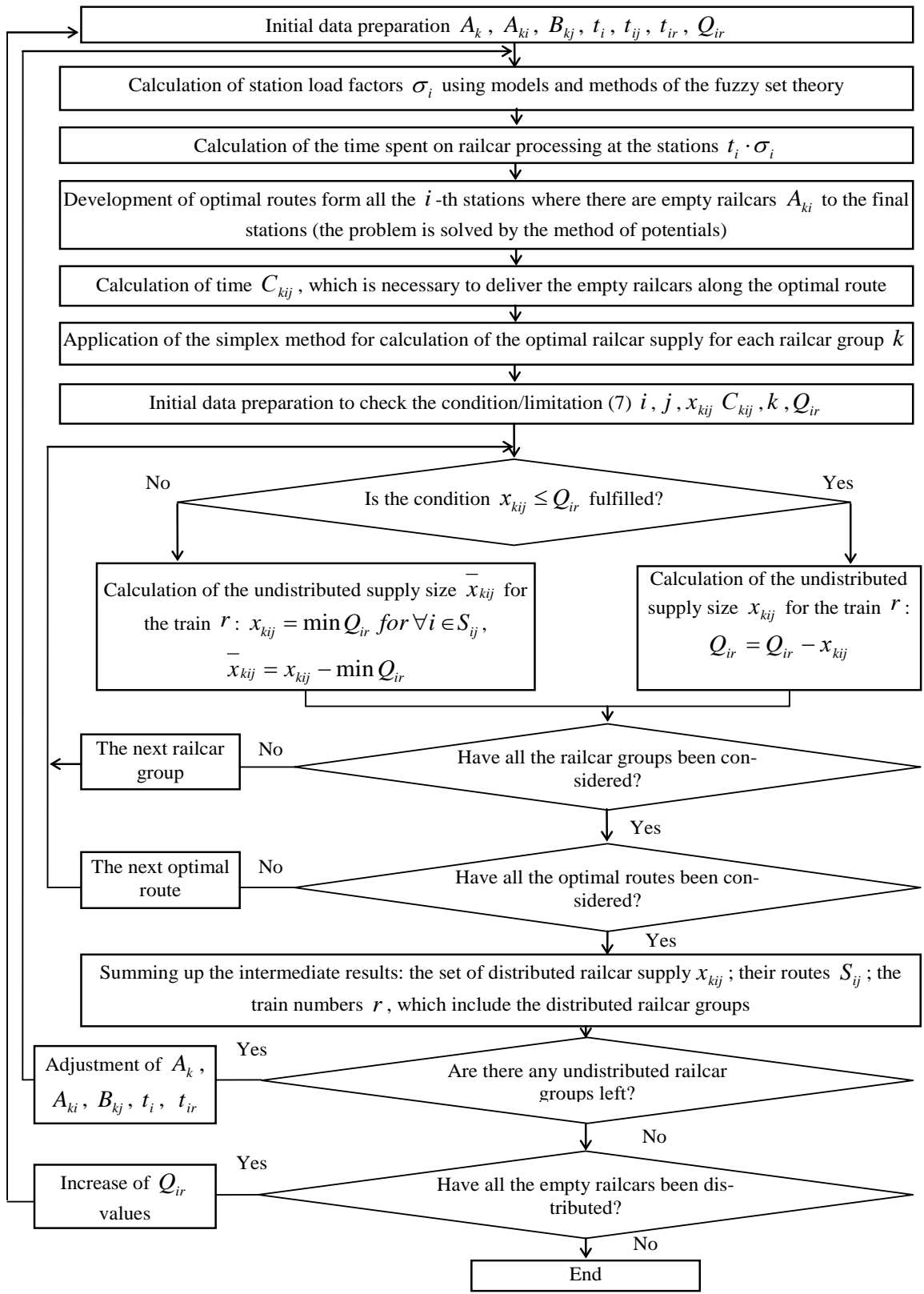

Fig. 1. Algorithm of optimal empty railcar distribution at railway transport nodes
Рис. 1. Алгоритм оптимального распределения порожних вагонов в железнодорожных транспортных
узлах

7. If the condition (7) is not fulfilled for the train r on any of the route vertex, the size of railcar supply $x_{kij}$ is considered equal to the minimum value of $x_{kij} = \min Q_{ir}$ for $\forall i \in S_{ij}$, and the difference $\overline{x}_{kij} = x_{kij} - \min Q_{ir}$ is stored as the undistributed supply.

If the condition (7) is fulfilled, the values of $Q_{ir}$ for all route vortexes are reduced by the railcar supply size $Q_{ir} = Q_{ir} - x_{kij}$ and the railcar supply $x_{kij}$ is stored as distributed. If after this the value of $Q_{ir}$ is equal to zero for the train $r$, this train is excluded from further calculations.

8. If there are several groups of railcars $k$ going along this route, the condition (7) is checked successively for each group. The order of consideration can be determined by priority in delivery of certain empty railcar groups if the enterprise has such a priority.

9. Calculation of the number of the undistributed railcars in each group as the sum of the number of railcars in the undistributed supply $A_{ki} = \sum_{i=1}^{M} \overline{x}_{kij}$ .

10. Summing up the intermediate results: the set of distributed railcar supply $x_{kij}$ ; their routes $S_{ij}$ ; the train numbers $r$ , which include the distributed railcar groups.

11. Repetition of the algorithm starting from the second stage until there is no undistributed railcar supplies left. If still there are such railcar supplies by the end of the base period, they will form the carry-over of the empty railcars for the following base period. This carryover can be eliminated if you increase the value of $Q_{ir}$ for the trains. After that it will be necessary to redevelop the plan of empty railcar distribution starting from the first step of the algorithm.

## 3. CHOICE OF SOLUTION METHOD AND RECOMMENDATIONS ON MODEL SOFTWARE IMPLEMENTATION

There are numerous specialized methods of solving the transport problem of linear programming, which makes it possible to reduce the number of iterations during hand calculations. Can be found several attempts of solving such problem of optimization for the main railway transport [11, 12]. However, taking into account the fact that this model will be implemented into the operating information system of the metallurgical enterprise and that the dimension of the problem is relatively small (within the day only about 100 loading areas require empty railcars delivery and the number of railcar groups does not exceed 50), the universal simplex method can be used in the developed model.

One of the program libraries of linear programming (API), for example, Linear Programming Library GIPALS32 [13], can be used for the automated model solution by the simplex method as well as for the integration of the model into the operating information system.

The research group believes that the algorithm of «The optimal routes table» described in the paper should be used to search for the least time route of the transport network and to calculate the time consumption $C_{kij}$ on delivery of empty railcars from the starting station of each route $i$ to the final $j$ -th stations.

## 4. CONCLUSION

The results of the developed model application are: set of values $x_{kij}$ determining the optimal number of railcars in each group, which must be delivered within the base period to a certain loading areas (stations); optimal routes of movement $S_{ij}$ , train numbers $r$ for each station, which include supplied empty railcars (assigning of empty railcar supply to the trains).

The implementation of the developed model resulted in significant labor intensity reduction of the on-line empty railcars distribution along the loading areas of the railway transport node and, as a result of the rational use of the empty railcar fleet, it is expected that the dead time of railcars on the industrial enterprise railway tracks will decrease by 15-20%.